\documentclass[10pt]{article}
\usepackage{amssymb,amsmath,amsthm}

\begin{document}
 \title{\bf Multiple periodic solutions for two classes of nonlinear difference systems involving classical $(\phi_1,\phi_2)$-Laplacian}
 \author {Xingyong Zhang \footnote{E-mail address:  zhangxingyong1@163.com},\ \ Liben Wang\\
       {\footnotesize   Department of  Mathematics, Faculty of Science, Kunming University of Science and Technology,}\\
        {\footnotesize  Kunming, Yunnan, 650500, P.R. China}\\}
 \date{}
 \maketitle{}

 \begin{center}
 \begin{minipage}{15cm}
 \par
 \small  {\bf Abstract:}  In this paper, we investigate the
 existence of multiple periodic solutions for two classes of nonlinear difference systems involving
 $(\phi_1,\phi_2)$-Laplacian. First, by using an important critical point theorem due to B. Ricceri, we establish an  existence theorem of three periodic solutions for
 the first nonlinear difference system with $(\phi_1,\phi_2)$-Laplacian and two parameters.  Moreover,  for the second nonlinear difference
 system with $(\phi_1,\phi_2)$-Laplacian, by  using the Clark's Theorem, we obtain  a
 multiplicity result of periodic solutions under a symmetric
 condition. Finally, two  examples are given to verify
 our theorems.

 \end{minipage}
 \end{center}

 \section{Introduction and main results}
 \setcounter{equation}{0}
  \allowdisplaybreaks

 \ \quad Let $\mathbb{R}$ denote the real numbers, $\mathbb{Z}$ the integers,
 Given $ a<b $ in $\mathbb{Z}$. Let
 $\mathbb{Z}[a,b]=\{a,a+1,....,b\}$. Let $ T>1 $ and $ N $ be fixed
 positive integers.

 Firstly, in this paper,  we are concerned with the existence of three periodic solutions
 for the following nonlinear difference system:
 \begin{eqnarray}\label{eq1}
 \left\{
  \begin{array}{l}
 \mu\Delta\Big[\rho_1(t-1)\phi_{1}\big(\Delta u_{1}(t-1)\big)\Big]-\mu\rho_3(t)\phi_{3}(u_1(t))+\nabla_{u_{1}} W\big(t,u_{1}(t),u_{2}(t)\big)=0
 \\
 \mu\Delta\Big[\rho_2(t-1)\phi_{2}\big(\Delta u_{2}(t-1)\big)\Big]-\mu\rho_4(t)\phi_{4}(u_2(t))+\nabla_{u_{2}} W\big(t,u_{1}(t),u_{2}(t)\big)=0,
    \end{array}
 \right.
 \end{eqnarray}
 where $\mu\in\mathbb R$, $\rho_i:\mathbb R\to \mathbb R^+$, $\phi_{i}$, $i=1,2,3,4$ satisfy the following conditions:
 \vskip2mm
 \noindent
 ($\mathcal{\rho}$) {\it $\rho_i$ are $T$-periodic and $\min_{t\in \mathbb Z[1,T]}\rho_i(t)>0$, $i=1,2,3,4$;}

 \vskip0mm
 \noindent
 ($\mathcal{A}1$) {\it $\phi_i:\mathbb R^N\to \mathbb R^N$ are
 homeomorphisms such that $\phi_i(0)=0$, $\phi_i=\nabla \Phi_i$,  with
 $\Phi_i\in C^1(\mathbb R^N,[0,+\infty))$ strictly
 convex and $\Phi_i(0)=0$, $i=1,2,3,4$.}
 \vskip2mm
 \noindent
 {\bf Remark} Assumption ($\mathcal{A}1$) is given in \cite{Mawhin2012}, where it is used to characterize the classical homeomorphism.

 \vskip2mm
 \noindent
 Moreover, assume that
 \vskip2mm
 \noindent
 ($\mathcal{A}2$) {\it  $W(t,x_1,x_2)=F(t,x_1,x_2)-\lambda G(t,x_1,x_2)+\nu H(t,x_1,x_2)$,
 where $\lambda, \nu\in\mathbb R$,
 $F,G,H:\mathbb{Z}\times \mathbb{R}^N\times \mathbb{R}^N\longrightarrow \mathbb{R}^N, (t,x_{1},x_{2})\longrightarrow F(t,x_{1},x_{2})$, $(t,x_{1},x_{2})\longrightarrow G(t,x_{1},x_{2})$,  $(t,x_{1},x_{2})\longrightarrow H(t,x_{1},x_{2})$ are $T$-periodic in
 $t$ for all $(x_{1},x_{2})\in \mathbb{R}^N\times\mathbb{R}^N $ and
 continuously differentiable in $(x_{1},x_{2})$ for every
 $t\in\mathbb{Z}[1,T]$.}
 \par
 It is well known that variational methods have been  important tools to study the existence and multiplicity of solutions for various difference
 systems. Lots of contributions has been obtained (for example, see
 \cite{Mawhin2012}, \cite{BC2009}-\cite{ZTZ}, \cite{TL}). Recently, in \cite{Mawhin2012} and \cite{Mawhin2013}, by using a variational approach,  Mawhin
 investigated the following second order nonlinear difference
 systems with $\phi$-Laplacian:
 \begin{eqnarray}
 \label{a1}  \Delta \phi[\Delta u(n-1)]=\nabla_u F[n,u(n)]+h(n) \quad (n\in \mathbb  Z),
 \end{eqnarray}
 where $\phi=\nabla \Phi$, $\Phi$ strictly convex, is a homeomorphism
 of $\mathbb R^N$ onto the ball $B_a\subset\mathbb R^N$ or of $B_a$
 onto $\mathbb R^N$.
 The assumption about $\phi$ implies three cases: firstly, classical homeomorphism if $\phi:\mathbb R^N\to\mathbb R^N$, for
 example, $\phi(0)=0$, $\phi(x)=|x|^{p-2}x$ for some $p>1$ and all $x\in \mathbb R^N/\{0\}$; secondly, bounded
 homeomorphism if $\phi:\mathbb R^N\to B_a$ $(a<+\infty)$, for example, $\phi(x)=\frac{x}{\sqrt{1+|x|^2}}\in B_1$ for all $x\in \mathbb
 R^N$; finally, singular
 homeomorphism if $\phi: B_a\subset\mathbb R^N\to\mathbb R^N$, for example, $\phi(x)=\frac{x}{\sqrt{1-|x|^2}}$ for all $x\in B_1$.   Under some reasonable assumptions, by using variational approach,
 Mawhin obtained system (\ref{a1}) has at least one
 $T$-periodic solution or $N+1$ geometrically distinct $T$-periodic solutions.
 \par
 However, to the best of our knowledge, except for recent works in \cite{WZ2014} and \cite{ZW2015} which are made by our first author and his cooperator named Yun Wang, there are no
 people to investigate the existence and multiplicity of   solutions  for system  involving classical $(\phi_1,\phi_2)$-Laplacian. In \cite{WZ2014}, Wang and our first author investigated the multiplicity of $T$-periodic solutions for the following nonlinear difference system:
  \begin{eqnarray}\label{eeq1}
 \left\{
  \begin{array}{l}
 \Delta \phi_{1}\big(\Delta u_{1}(t-1)\big)=\nabla_{u_{1}} F\big(t,u_{1}(t),u_{2}(t)\big)+h_{1}(t)\\
 \Delta \phi_{2}\big(\Delta u_{2}(t-1)\big)=\nabla_{u_{2}} F\big(t,u_{1}(t),u_{2}(t)\big)+h_{2}(t),
    \end{array}
 \right.
 \end{eqnarray}
 where $F:{\mathbb Z}\times {\mathbb R}^{N}\times{\mathbb R}^{N}\to {\mathbb R}$ and $\phi_{m}, m=1,2$ satisfy the following
 condition:
 \vskip1mm
 \noindent
 ($\mathcal{A}$)\ \  $\phi_{i}$ is a homeomrphism from ${\mathbb R}^{N}$ onto $B_{a}\subset {\mathbb R}^{N}(a\in(0,+\infty])$, such that
       $
       \phi_{i}(0)=0, \phi_{i}=\nabla\Phi_{i},
       $
       with $ \Phi_{i}\in C^{1}({\mathbb R}^{N},[0,+\infty])$ strictly convex and
       $\Phi_{i}(0)=0$, $m=1,2$.
 \vskip2mm
 \noindent
 Assumption ($\mathcal{A}$) implies that $\Phi_{i}$, $i=1,2$ are the classical
 homeomorphisms or the bounded homeomorphisms. They  investigated the case that $F(t,x_1,x_2)$ is periodic on $r_1$ components of variables
 $x_{1}^{(1)},\cdots,x_{N}^{(1)}$ and $r_2$ components of variables
 $x_{1}^{(2)},\cdots,x_{N}^{(2)}$, where $1\le r_1\le N$ and $1\le r_2\le N$.
  By using a critical point theorem in \cite{MW1989} and a generalized saddle point theorem in \cite{Liu}, they  obtain that system (\ref{eeq1}) has at least $r_1+r_2+1$  geometrically distinct $T$-periodic solutions.
  Their results generalize those corresponding to classical homeomorphism and  bounded homeomorphism in \cite{Mawhin2013}.
  \par
 In \cite{ZW2015}, our first author and Wang investigated  the existence  of  homoclinic solutions for the following nonlinear difference  systems involving classical ($\phi_1$, $\phi_2$)-Laplacian:
  \begin{eqnarray}\label{eeeq1}
 \left\{
  \begin{array}{l}
 \Delta \phi_{1}\big(\Delta u_{1}(t-1)\big)+\nabla_{u_{1}} V\big(t,u_{1}(t),u_{2}(t)\big)=f_{1}(t)\\
 \Delta \phi_{2}\big(\Delta u_{2}(t-1)\big)+\nabla_{u_{2}} V\big(t,u_{1}(t),u_{2}(t)\big)=f_{2}(t),
    \end{array}
 \right.
 \end{eqnarray}
 where $t\in\mathbb Z$, $u_m(t)\in \mathbb R^N$, $m=1,2$, $V(t,x_1,x_2)=-K(t,x_1,x_2)+W(t,x_1,x_2)$, $K,W:\mathbb{Z}\times\mathbb{R}^{N}\times\mathbb{R}^{N}\to \mathbb R$ and $\phi_{m}, m=1,2$ satisfy
 assumption  ($\mathcal{A}1$). They first  improve some inequalities in \cite{HC2011}.  Then by using a linking theorem in \cite{Schechter},
 some new existence results of homoclinic solutions  for system (\ref{eeeq1}) are obtained when $W$ has super $p$-linear growth and $K$ has sub $p$-linear growth.
 \par
   Inspired by \cite{Mawhin2012},  \cite{Ricceri}, \cite{WZ2014}, \cite{ZW2015} and \cite{Mawhin2013},  in this paper, we are interested in
 the existence of three $T$-periodic solutions for system
 (\ref{eq1}).  By using an
 important three critial point theorem established by B. Ricceri in \cite{Ricceri}, we
 investigate the existence of three $T$-periodic solutions for
 system (\ref{eq1}),
  as stated
 in the following.
 \par
 Define
   \begin{eqnarray*}
  &   &  I(u)   =    \sum\limits_{t = 1}^{T}\left[\rho_1(t)\Phi_1(\Delta u_{1}(t))+\rho_2(t)\Phi_2(\Delta u_{2}(t))+\rho_3(t)\Phi_3(u_1(t))+\rho_4(t)\Phi_4(u_2(t))\right],\\
  &   &  \Psi(u)   =   -  \sum\limits_{t = 1}^{T}  F(t,u_{1}(t),u_{2}(t)),\quad \Phi(u)   =     \sum\limits_{t = 1}^{T}  G(t,u_{1}(t),u_{2}(t)),\\
   &   &   \Gamma(u)   =   - \sum\limits_{t = 1}^{T}  H(t,u_{1}(t),u_{2}(t)), \quad\quad u\in   E,
 \end{eqnarray*}
  where the definitions of $E$ and its norm are in section 2 below.
  \vskip1mm
 \noindent
 {\bf Theorem 1.1.} {\it Suppose that ($\mathcal{\rho}$),  ($\mathcal{A}1$),  ($\mathcal{A}2$) and the following conditions hold:\\
  ($\mathcal{A}3$) there exist positive constants $c_i$
  $(i=1,2,3,4)$, $\theta>1$  such that
  \begin{eqnarray*}
    (\phi_i(x)-\phi_i(y),x-y)\ge c_i|x-y|^\theta,\forall\  x,y\in \mathbb    R^N,i=1,2,3,4,
  \end{eqnarray*}
  where $(\cdot,\cdot)$ stands for the usual product in $\mathbb  R^N$;\\
  ($\mathcal{A}4$) $\lim_{|x|\to\infty}\Phi_i(x)=+\infty$ and  there exist positive constants $l\ge \theta$, $d_i$ and $m_i$ such that
   $\Phi_i(x)\le d_i|x|^l+m_i$ for all $x\in \mathbb R^N$, $(i=1,2,3,4)$;\\
  ($\mathcal{A}5$)  for all $t\in \mathbb  Z[1,T]$ and all $\lambda>0$, there exists $C_0(\lambda)\in \mathbb R$  such that  for  all $(x_1,x_2)\in \mathbb R^N\times \mathbb
  R^N$,
  \begin{eqnarray*}
  \lim_{|x_1|+|x_2|\to \infty}\frac{F(t,x_1,x_2)}{|x_1|^l+|x_2|^l}=+\infty,\quad \lambda G(t,x_1,x_2)\ge  F(t,x_1,x_2)+C_0(\lambda);
  \end{eqnarray*}
  ($\mathcal{A}6$)$\sum_{t=1}^T G(t,0,0)= 0$.\\
 Then  for each $r>0$, for each $\mu>\max\{0,\mu^*(I,\Psi,\Phi,r)\}$, and
 for each compact interval $[a,b]\subset]0,\beta(\mu
 I+\Psi,\Phi,r)[$, there exists a number $\rho>0$ with the following
 property: for every $\lambda\in [a,b]$, there exists
 $\delta>0$ such that, for each $\nu\in [0,\delta]$, system (1.1)
 has at least three $T$-periodic solutions in $E$ whose norms are less than
 $\rho$, where }
 \begin{eqnarray*}
 &   &  \beta(\mu I+\Psi,\Phi,r)=\sup_{u\in\Phi^{-1}(]r,+\infty[)}\frac{\mu  I(u)+\Psi(u)-\inf\limits_{\Phi^{-1}(]-\infty,r])}(\mu I+\Psi)}{r-\Phi(u)}\\
 &   &  \mu^*(I,\Psi,\Phi,r)=\inf\left\{\frac{\Psi(u)-\gamma+r}{\eta_r-I(u)}:u\in E,\Phi(u)<r, I(u)<\eta_r\right\}\\
 &   &  \gamma=\inf\limits_E(\Psi(u)+\Phi(u)),\quad \eta_r=\inf\limits_{u\in \Phi^{-1}(r)}I(u).
 \end{eqnarray*}

 \vskip2mm
 Inspired by \cite{Ricceri}, we have the following corollary:
 \vskip2mm
 \noindent
 {\bf Corollary 1.1.} {\it Suppose that  ($\mathcal{\rho}$), ($\mathcal{A}1$)-($\mathcal{A}4$) and ($\mathcal{A}6$) hold. If\\
 ($\mathcal{A}5$)$'$ there exists  $s>l$ such that for every $t\in \mathbb Z[1,T]$,
  \begin{eqnarray*}
   \lim_{|x_1|+|x_2|\to
   \infty}\frac{F(t,x_1,x_2)}{|x_1|^l+|x_2|^l}=+\infty,\quad  \lim_{|x_1|+|x_2|\to
   \infty}\frac{F(t,x_1,x_2)}{|x_1|^s+|x_2|^s}<+\infty
  \end{eqnarray*}
  and
  $$
   \quad  \lim_{|x_1|+|x_2|\to
   \infty}\frac{G(t,x_1,x_2)}{|x_1|^s+|x_2|^s}=+\infty,
  $$
 then the conclusion of Theorem 1.1 holds.
 }
 \vskip2mm
 \noindent
 {\bf Remark 1.1.} {\it There exist  examples satisfying ($\mathcal{A}1$)-($\mathcal{A}6$) in Theorem 1.1. For example, let $T>1$ and $N$ be fixed integer. Let $\theta\ge 2$ and $q_i\ge 2$, $i=1,2,3,4$.  Assume that $\phi_1(y)=|y|^{\theta-2}y+|y|^{q_1-2}y$,
 $\phi_2(y)=|y|^{\theta-2}y+ |y|^{q_2-2}y$,
 $\phi_3(y)=|y|^{\theta-2}y+ |y|^{q_3-2}y$,
 $\phi_4(y)=|y|^{\theta-2}y+ |y|^{q_4-2}y$,
 $\rho_i$ are $T$-periodic and satisfy
 $\rho_i> 0$ for all $t\in \mathbb Z[1,T]$, $i=1,2,3,4$.
 Then  $\Phi_1(y)=\frac{|y|^{\theta}}{\theta}+ \frac{|y|^{q_1}}{q_1}$,  $\Phi_2(y)=\frac{|y|^{\theta}}{\theta}+
 \frac{|y|^{q_2}}{q_2}$, $\Phi_3(y)=\frac{|y|^{\theta}}{\theta}+
 \frac{|y|^{q_3}}{q_3}$,  $\Phi_4(y)=\frac{|y|^{\theta}}{\theta}+
 \frac{|y|^{q_4}}{q_4}$.
 \par
 Note that
 $$
 (|x|^{\theta-2}x-|y|^{\theta-2}y,x-y)\ge c|x-y|^\theta
 $$
 for all $x,y\in\mathbb R^N$, $\theta\ge 2$ and some $c>0$ (see \cite{simon}).
 Hence,
 \begin{eqnarray*}
 &    & (\phi_1(x)-\phi_1(y),x-y)\\
 &  = &  (|x|^{\theta-2}x+|x|^{q_1-2}x-|y|^{\theta-2}y- |y|^{q_1-2}y,x-y)\\
 &  = &   (|x|^{\theta-2}x-|y|^{\theta-2}y,x-y) +( |x|^{q_1-2}x- |y|^{q_1-2}y,x-y)\\
 &\ge &  (|x|^{\theta-2}x-|y|^{\theta-2}y,x-y)\\
 &\ge &  c |x-y|^\theta,\quad \forall\ x,y\in \mathbb R^N.
  \end{eqnarray*}
 Similarly, we have
 \begin{eqnarray*}
 (\phi_i(x)-\phi_i(y),x-y)\ge c_i |x-y|^\theta , \quad \forall\ x,y\in \mathbb
 R^N, i=2,3,4
 \end{eqnarray*}
 for some $c_i>0$, $i=2,3,4$. So $(\mathcal{A}3)$ holds.
 \par
 Take
 $l=\max\{\theta,q_1,q_2,q_3,q_4\}$ and let
  \begin{eqnarray}
  \label{ddd1}         F(t,x_1,x_2)
 &  =  &  \left[(\cos^2\frac{\pi t}{T}+2)|x_1|^{l}+(|\cos\frac{\pi t}{T}|+2)|x_2|^{l}\right]\ln(|x_1|^2+|x_2|^2+1)\nonumber\\
 G (t,x_1,x_2)& =  & (|\sin\frac{\pi t}{T}|+2)(|x_1|^{l}+|x_2|^{l})^2\ln(|x_1|^2+|x_2|^2+1)\nonumber\\
 H (t,x_1,x_2)& =  & (\cos^2\frac{\pi t}{T}+2)\sin(|x_1|^2+|x_2|^2+2)\nonumber\\
 W(t,x_1,x_2) & =  & F(t,x_1,x_2)-\lambda G (t,x_1,x_2)+\nu H (t,x_1,x_2).\nonumber
 \end{eqnarray}
 Then it is easy to obtain that ($\mathcal{A}2$) and ($\mathcal{A}6$) hold and $\Phi_i$ satisfy ($\mathcal{A}1$) and ($\mathcal{A}4$),
 $i=1,2,3,4$. Moreover,
  \begin{eqnarray*}
             \lim_{|x_1|+|x_2|\to\infty} \frac{F (t,x_1,x_2)}{|x_1|^l+|x_2|^l}
     \ge  2\lim_{|x_1|+|x_2|\to\infty} \ln(|x_1|^2+|x_2|^2+1)
      =   +\infty
  \end{eqnarray*}
 and for all $\lambda>0$,
  \begin{eqnarray*}
             \lim_{|x_1|+|x_2|\to\infty} \frac{\lambda G (t,x_1,x_2)}{F(t,x_1,x_2)}
    &  =  & \lim_{|x_1|+|x_2|\to\infty} \frac{\lambda \left(|\sin\frac{\pi t}{T}|+2\right)(|x_1|^{l}+|x_2|^{l})^2}{\left[(\cos^2\frac{\pi t}{T}+2)|x_1|^{l}+(|\cos\frac{\pi t}{T}|+2)|x_2|^{l}\right]}\\
    & \ge & \lim_{|x_1|+|x_2|\to\infty}    \frac{2\lambda(|x_1|^{l}+|x_2|^{l})^2}{3\left(|x_1|^{l}+|x_2|^{l}\right)}\\
    &  =  & \lim_{|x_1|+|x_2|\to\infty}    \frac{2\lambda}{3}\left(|x_1|^{l}+|x_2|^{l}\right)=+\infty.
  \end{eqnarray*}
  Hence, ($\mathcal{A}5$) holds.
}

 Moreover, in this paper, we are also concerned with the  multiplicity of  $T$-periodic solutions
 for the following nonlinear difference system:
 \begin{eqnarray}\label{eqq1}
 \left\{
  \begin{array}{l}
 \Delta\left(\gamma_1(t-1) \phi_{1}\big(\Delta u_{1}(t-1)\big)\right)-\gamma_3(t)\phi_{3}(|u_1(t)|)+\nabla_{u_{1}} F\big(t,u_{1}(t),u_{2}(t)\big)=0
 \\
 \Delta\left(\gamma_2(t-1) \phi_{2}\big(\Delta u_{2}(t-1)\big)\right)-\gamma_4(t)\phi_{4}(|u_2(t)|)+\nabla_{u_{2}} F\big(t,u_{1}(t),u_{2}(t)\big)=0,
    \end{array}
 \right.
 \end{eqnarray}
  where  $\gamma_i:\mathbb R\to \mathbb R^+$ satisfy the following conditions:
 \vskip2mm
 \noindent
 ($\mathcal{\gamma}$) {\it $\gamma_i$ are $T$-periodic and $\min_{t\in \mathbb Z[1,T]}\gamma_i(t)>0$, $i=1,2,3,4$,}
 \vskip2mm
 \noindent
  and $\phi_{i}$, $i=1,2,3,4$ satisfy the assumption ($\mathcal{A}1$)
  and the following condition:
 \vskip2mm
 \noindent
 ($\mathcal{\phi}$) {\it there exist positive constants $p>1,q>1$, $a_i$, $b_i$, $i=1,2,3,4$ such that
 $$
  a_i|x|^q\le \Phi_i(x)\le b_i |x|^q, i=1,3,\quad \forall\  x\in  \mathbb R^N
 $$
 and}
 $$
  a_i|x|^p\le \Phi_i(x)\le b_i |x|^p, i=2,4, \quad \forall\  x\in   \mathbb  R^N.
 $$
  \vskip2mm
 \noindent
 Moreover, $F:\mathbb{Z}\times \mathbb{R}^N\times \mathbb{R}^N\longrightarrow
 \mathbb{R}, (t,x_{1},x_{2})\longrightarrow F(t,x_{1},x_{2})$ is
 $T$-periodic in
 $t$ for all $(x_{1},x_{2})\in \mathbb{R}^N\times\mathbb{R}^N $ and
 continuously differentiable in $(x_{1},x_{2})$ for every
 $t\in\mathbb{Z}[1,T]$.
 \par
 When $\Phi_i(x)=\frac{1}{q}|x|^q$, $i=1,3$ and  $\Phi_i(x)=\frac{1}{p}|x|^p$,
 $i=2,4$,
  system (\ref{eqq1}) can be seen as a discrete
 analogue of the following $(q,p)$-Laplacian differential systems:
 \begin{eqnarray}\label{eqq2}
 \left\{
  \begin{array}{l}
 \dfrac{d\left(\gamma_1(t) |\dot{u}_{1}(t)|^{q-2}\dot{u}_{1}(t)\right)}{dt}-\gamma_3(t)|u_1(t)|^{q-2}u_1(t)+\nabla_{u_{1}} F\big(t,u_{1}(t),u_{2}(t)\big)=0
 \\
 \dfrac{d\left(\gamma_2(t) |\dot u_{2}(t)|^{p-2}\dot u_{2}(t)\right)}{dt}-\gamma_4(t)|u_2(t)|^{p-2}u_2(t)+\nabla_{u_{2}} F\big(t,u_{1}(t),u_{2}(t)\big)=0.
    \end{array}
 \right.
 \end{eqnarray}
 Recently, by using variational methods, system (\ref{eqq2}) has been investigated by some authors
 (for example, see \cite{YC2012}-\cite{LZ2011}) and some interesting results
 on the existence and multiplicity of solutions have been obtained. However, to the best of our knowledge, there are no
 people to investigate the nonlinear difference  system (\ref{eqq1}). In this paper, inspired by
 \cite{TL}-\cite{LZ2011}, we are interested in
 the existence and multiplicity of $T$-periodic solutions for system
 (\ref{eqq1}). By using the  Clark's
 theorem, we obtain the following theorem.
  \vskip2mm
 \noindent
 {\bf Theorem 1.2.} {\it Suppose that ($\mathcal{\gamma}$), ($\mathcal{\phi}$) and  the following conditions hold:\\
 ($\mathcal{F}$0) there exist $\alpha_1\in
 [0,q)$, $\alpha_2\in [0,p)$, $h_i:\mathbb Z[1,T]\to \mathbb R^+, i=1,2$ and $l:\mathbb Z[1,T]\to \mathbb R^+$ such that
  $$
    F(t,x_1,x_2)\le h_1(t) |x_1|^{\alpha_1}+h_2(t) |x_2|^{\alpha_2}+l(t).
  $$
  ($\mathcal{F}$1) $F(t,0,0)=0$;\\
  ($\mathcal{F}$2) $F(t,-x_1,-x_2)=F(t,x_1,x_2)$;\\
  ($\mathcal{F}$3) there exist constants $\beta_i\in (1,\min\{q,p\})$, $M_i\in (0,\infty)$,  $i=1,2$ and $\delta\in (0,1)$ such that
   $$
     F(t,x_1,x_2)\ge M_1|x_1|^{\beta_1}+ M_2|x_2|^{\beta_2},\quad \forall \  |x_1|<\delta,|x_2|<\delta.
   $$
  Then system (\ref{eqq1}) has at least $2NT$ distinct pairs of nonzero
 solutions.}

 \vskip-2mm
 \noindent
 \section{Preliminaries}
 \setcounter{equation}{0}
 \quad At first, we make some preliminaries.
 Define
 \begin{eqnarray*}
 E_T=\{h:=\{h(t)\}_{t\in \mathbb Z} |h(t+T)=h(t),h(t)\in\mathbb{R}^N, t\in\mathbb{Z}\}
 \end{eqnarray*}
  and let $E=E_T\times E_T$.
 For $h\in E_T$, set
 \begin{eqnarray}
 \label{ffffff1}\|h\|_{r}=\left(\sum\limits_{t = 1}^{T}|h(t)|^{r}\right)^{1/r} \ \
 \mbox{and}\ \
 \|h\|_{\infty}=\max_{t\in\mathbb{Z}[1,T]}|h(t)|,\quad
 r>1.
 \end{eqnarray}
 Obviously, we have
 \begin{eqnarray}
 \label{ff1} \|h\|_{\infty}\le \|h\|_{r}\le T^{1/r}  \|h\|_{\infty}.
 \end{eqnarray}
 On $E_T$, we define
  \begin{eqnarray*}
     \|h\|_{E_T}=\left(\sum_{t=1}^T|\Delta  h(t)|^\theta+\sum_{t=1}^T|h(t)|^\theta\right)^{1/\theta}
  \end{eqnarray*}
  and
  \begin{eqnarray*}
     \|h\|_{[E_T]}=\left(\sum_{t=1}^T|\Delta   h(t)|^l+\sum_{t=1}^T|h(t)|^l\right)^{1/l}
  \end{eqnarray*}
 For $u=(u_1,u_2)\in E$, define
 $$
 \|u\|=\|u_1\|_{E_T}+\|u_2\|_{E_T}.
 $$
 Then $E$ is a separable  and reflexive Banach space.
  Moreover,   $\|\cdot\|_{E_T}$  is equivalent to
  $\|\cdot\|_{r} (r>1)$ and $\|\cdot\|_{[E_T]}$. Hence, there exist positive constants $C_i$ $(i=1,\cdots, 6)$ such that
 \begin{eqnarray}
        \label{d2}
  &   &  C_1\|\cdot\|_{E_T}\le \|\cdot\|_{\theta}\le C_2\|\cdot\|_{E_T},\\
   \label{d3} &   &  C_3\|\cdot\|_{E_T}\le \|\cdot\|_{l}\le C_4\|\cdot\|_{E_T},\\
 \label{d4}  &   &  C_5\|\cdot\|_{E_T}\le \|\cdot\|_{[E_T]}\le C_{6}\|\cdot\|_{E_T}.
 \end{eqnarray}

 \vskip2mm
 \noindent
 {\bf Lemma 2.1 }(see \cite{WZ2014}) {\it Let $L:\mathbb{Z}[1,T]\times\mathbb{R}^{N}\times\mathbb{R}^{N}\times\mathbb{R}^{N}\times\mathbb{R}^{N}\longrightarrow\mathbb{R},
 (t,x_{1},x_{2},y_{1},y_{2})\longrightarrow L(t,x_{1},x_{2},y_{1},y_{2})$  and assume that $L$ is continuously differentiable in
 $(x_{1},x_{2},y_{1},y_{2})$ for all $ t\in\mathbb{Z}[1,T]$.  Then the functional $\varphi:E\to \mathbb{R} $ defined by
 $$
   \varphi(u)= \varphi(u_{1},u_{2})=\sum\limits_{t=1}^{T}L\big(t,u_{1}(t),u_{2}(t),\Delta u_1(t),\Delta u_2(t)\big)
 $$
 is continuously differentiable on $E$ and for $u,v\in E$,}
 \begin{eqnarray*}
  &     &    \langle\varphi'(u),v\rangle
    =    \langle\varphi'(u_{1},u_{2}),(v_{1},v_{2})\rangle\\
  &  =  &  \sum\limits_{t=1}^{T}\left[\big(D_{x_{1}}L(t,u_{1}(t),u_{2}(t),\Delta u_{1}(t),\Delta  u_{2}(t)),v_{1}(t)\big)\right.\\
  &     & ¡¡+\big(D_{y_{1}}L(t,u_{1}(t),u_{2}(t),\Delta u_{1}(t),\Delta u_{2}(t)),\Delta v_{1}(t)\big)\\
  &     &  +\big(D_{x_{2}}L(t,u_{1}(t),u_{2}(t),\Delta u_{1}(t),\Delta  u_{2}(t)),v_{2}(t)\big)\\
  &     &  \left.+\big(D_{y_{2}}L(t,u_{1}(t),u_{2}(t),\Delta u_{1}(t),\Delta u_{2}(t)),\Delta  v_{2}(t)\big)\right].
  \end{eqnarray*}

 \par
  Let
  \begin{eqnarray*}
 L(t,x_{1},x_{2},y_{1},y_{2})&  =  & \mu[\rho_1(t)\Phi_1(y_{1})+\rho_2(t)\Phi_2(y_{2})+\rho_3(t)\Phi_3(x_{1})+\rho_4(t)\Phi_4(x_{2})]\\
                             &     &   - F(t,x_{1},x_{2})+\lambda  G(t,x_{1},x_{2})-\nu H(t,x_{1},x_{2}),
 \end{eqnarray*}
 where $ F,G,H $: $\mathbb{Z}[1,T]\times\mathbb{R}^N\times\mathbb{R}^N\longrightarrow\mathbb{R}^N$ are continuously differentiable in
 $(x_{1},x_{2})\in \mathbb{R}^N\times\mathbb{R}^N$ for all $t\in\mathbb{Z}[1,T]$.
 Then
 \begin{eqnarray}
 \label{b1}         \varphi(u)
    &  =  &     \sum\limits_{t = 1}^{T}\left[\mu(\rho_1(t)\Phi_1(\Delta u_1(t))+\rho_2(t)\Phi_2(\Delta u_2(t))+\rho_3(t)\Phi_3(u_1(t))+\rho_4(t)\Phi_4(u_2(t)))\right.\nonumber\\
    &     &     \left.-F(t,u_{1}(t),u_{2}(t))+\lambda G(t,u_{1}(t),u_{2}(t))-\nu H(t,u_{1}(t),u_{2}(t))\right].
 \end{eqnarray}
 Obviously, when
 ($\mathcal{A}1$) and ($\mathcal{A}2$) hold, $\varphi$ is continuously differentiable on $E$
 and for $ \forall\ u,v\in E$, we have
 \begin{eqnarray}
  \label{ff2}&     & \langle\varphi'(u),v\rangle
              =    \langle\varphi'(u_{1},u_{2}),(v_{1},v_{2})\rangle\nonumber\\
            &  =  &  \sum\limits_{t = 1}^{T}\left[\mu\rho_1(t)(\phi_1(\Delta u_1(t)),\Delta v_1(t))+\mu\rho_2(t)(\phi_2(\Delta u_2(t)),\Delta v_2(t))\right.\nonumber\\
            &     &  \left.+\mu\rho_3(t)(\phi_3( u_1(t)), v_1(t))+\mu \rho_4(t)(\phi_4( u_2(t)), v_2(t))]\right.\nonumber\\
            &     &  -\sum\limits_{t=1}^{T}\left[(\nabla_{u_1} F(t,u_{1}(t),u_{2}(t)),v_1(t))+(\nabla_{u_2}F(t,u_{1}(t),u_{2}(t)),v_2(t))\right]\nonumber\\
            &     &  +\lambda\sum\limits_{t=1}^{T}\left[(\nabla_{u_1} G(t,u_{1}(t),u_{2}(t)),v_1(t))+(\nabla_{u_2}G(t,u_{1}(t),u_{2}(t)),v_2(t))\right]\nonumber\\
             &     &  -\nu\sum\limits_{t=1}^{T}\left[(\nabla_{u_1} H(t,u_{1}(t),u_{2}(t)),v_1(t))+(\nabla_{u_2}H(t,u_{1}(t),u_{2}(t)),v_2(t))\right].
  \end{eqnarray}

 \noindent
 {\bf Lemma 2.2.} {\it If $u\in E$ is a solution of Euler equation
 $\varphi'(u)=0$, then  $u$ is a solution of system (1.1).}
 \vskip2mm
 \noindent
 {\bf Proof } At first, for any  $u=(u_1,u_2),v=(v_1,v_2)\in E$, we can obtain the following two   equalities:
  \begin{eqnarray}
 \label{b2}&  &  -\sum_{t=1}^T\left(\Delta\Big[\rho_1(t-1)\phi_1(\Delta u_1(t-1))\Big], v_1(t)\right) = \sum_{t=1}^T  (\rho_1(t)\phi_1(\Delta u_1(t)),\Delta v_1(t)),\\
 \label{b3}&  &   -\sum_{t=1}^T\left(\Delta\Big[\rho_2(t-1)\phi_2(\Delta u_2(t-1))\Big], v_2(t)\right) = \sum_{t=1}^T  (\rho_2(t)\phi_2(\Delta u_2(t)),\Delta v_2(t)).
  \end{eqnarray}
  In fact, since $u_1(t)=u_1(t+T)$ and $v_1(t)=v_1(t+T)$ for all $t\in \mathbb  Z$, then
 \begin{eqnarray*}
 &     &   -\sum_{t=1}^T\left(\Delta\Big[\rho_1(t-1)\phi_1(\Delta u_1(t-1))\Big], v_1(t)\right) \\
 &  =  &   -\sum_{t=1}^{T}(\rho_1(t)\phi_1(\Delta u_1(t)), v_1(t))+\sum_{t=1}^T(\rho_1(t-1)\phi_1(\Delta u_1(t-1)), v_1(t)) \\
 &  =  &   -\sum_{t=1}^{T}(\rho_1(t)\phi_1(\Delta u_1(t)), v_1(t))+\sum_{t=1}^{T-1}(\rho_1(t)\phi_1(\Delta u_1(t)), v_1(t+1))+ (\rho_1(0)\phi_1(\Delta u_1(0)), v_1(1))\\
 &  =  &   \sum_{t=1}^{T}(\rho_1(t)\phi_1(\Delta u_1(t)), \Delta v_1(t))+ (\rho_1(0)\phi_1(\Delta u_1(0)), v_1(1))-(\rho_1(T)\phi_1(\Delta u_1(T)), v_1(T+1))\\
 &  =  &   \sum_{t=1}^{T}(\rho_1(t)\phi_1(\Delta u_1(t)), \Delta v_1(t)).
 \end{eqnarray*}
 Hence, (\ref{b2}) holds. Similarly, it is easy to get (\ref{b3}). Since $\varphi'(u)=0$,
 then for all $v=(v_1,0)\in E$, (\ref{ff2}) implies that
 \begin{eqnarray}
  \label{b5}&     &       \sum\limits_{t = 1}^{T}\left[\mu(\rho_1(t)\phi_1(\Delta u_1(t)),\Delta  v_1(t))+\mu(\rho_3(t)\phi_3( u_1(t)), v_1(t))\right]\nonumber\\
            &  =  &  \sum\limits_{t=1}^{T} (\nabla_{u_1} F(t,u_{1}(t),u_{2}(t)),v_1(t))
                     -\lambda\sum\limits_{t=1}^{T} (\nabla_{u_1} G(t,u_{1}(t),u_{2}(t)),v_1(t))\nonumber\\
            &     &  +\nu \sum\limits_{t=1}^{T} (\nabla_{u_1} H(t,u_{1}(t),u_{2}(t)),v_1(t))
  \end{eqnarray}
 Note that $v_1$ is arbitrary. Then (\ref{b2}) and (\ref{b5}) imply that
 \begin{eqnarray*}
 \mu \Delta\Big[\rho_1(t-1)\phi_1(\Delta u_1(t-1))\Big]-\mu\rho_3(t)\phi_3( u_1(t))+\nabla_{u_{1}} W(t,u_{1}(t),u_{2}(t))=0.
 \end{eqnarray*}
 Similarly, Let $v_1=0$. We can obtain that
 \begin{eqnarray*}
  \mu \Delta\Big[\rho_2(t-1)\phi_2(\Delta u_2(t-1))\Big]-\mu\rho_4(t)\phi_4( u_2(t))+\nabla_{u_{2}} W(t,u_{1}(t),u_{2}(t))=0. \quad\Box
 \end{eqnarray*}

 \vskip2mm
 \par
 To prove Theorem 1.1, we will use the following three critical points theorem due to
 Ricceri \cite{Ricceri}.
 \vskip2mm
 \noindent
 {\bf Theorem 2.1} (see \cite{Ricceri})  {\it  Let $X$ be a reflexive real Banach space, $I:X\to \mathbb R$
 a sequentially weakly lower semicontinuous, coercive, bounded on
 each bounded subset of $X$, $C^1$ functional whose derivative admits a continuous inverse on $X^*$; $\Psi,\Phi:X\to\mathbb R$ two $C^1$ functionals with compact
 derivative. Assume also that the functional $\Psi+\lambda \Phi$ is bounded below for all $\lambda >0$ and that
  \begin{equation}
   \liminf_{\|x\|\to+\infty}\frac{\Psi(x)}{I(x)}=-\infty.
  \end{equation}
 Then, for each $r>\sup_M\Phi$, where $M$ is the set of all global
 minima of $I$, for each $\mu>\max\{0,\mu^*(I,\Psi,\Phi,r)\}$, and
 for each compact interval $[a,b]\subset]0,\beta(\mu
 I+\Psi,\Phi,r)[$, there exists a number $\rho>0$ with the following
 property: for every $\lambda\in [a,b]$ and every $C^1$ functional
 $\Gamma:X\to\mathbb R$ with compact derivative, there exists
 $\delta>0$ such that, for each $\nu\in [0,\delta]$, the equation
 $$
 \mu I'(x)+\Psi'(x)+\lambda \Phi'(x)+\nu \Gamma'(x)=0
 $$
 has at least three solutions in $X$ whose norms are less than
 $\rho$, where}
 \begin{eqnarray*}
 &   &  \beta(\mu I+\Psi,\Phi,r)=\sup_{x\in\Phi^{-1}(]r,+\infty[)}\frac{\mu  I(x)+\Psi(x)-\inf\limits_{\Phi^{-1}(]-\infty,r])}(\mu I+\Psi)}{r-\Phi(x)}\\
 &   &  \mu^*(I,\Psi,\Phi,r)=\inf\left\{\frac{\Psi(x)-\gamma+r}{\eta_r-I(x)}:x\in X,\Phi(x)<r,
 I(x)<\eta_r\right\}\\
 &   &  \gamma=\inf\limits_X(\Psi(x)+\Phi(x)),\quad \eta_r=\inf\limits_{x\in \Phi^{-1}(r)}I(x).
 \end{eqnarray*}

 \section{Proof of Theorem 1.1}
 \setcounter{equation}{0}
 \par
 \quad For the sake of convenience, we denote
 $$
  \rho_i^+=\max_{t\in \mathbb Z[1,T]}\rho_i(t),\quad \rho_i^-=\min_{t\in \mathbb Z[1,T]}  \rho_i(t), \ \ i=1,2,3,4.
 $$
 {\bf Proof of Theorem 1.1}\ \ We prove that $\varphi$ defined by (\ref{b1}) satisfies all the assumptions of Theorem 2.1. Let $X=E$. Then $E$ is a reflexive and separable Banach space.
 Since all  the topologies are equivalent in the finite dimensional Banach space $E$, then for any sequence $\{u^n\}\subset E$, assume that
  \begin{eqnarray}
  \label{fff1}u^n\to u^*  \mbox{ in } E \mbox{ as  } n\to\infty,
  \end{eqnarray}
 that is,
 \begin{eqnarray*}
  &     &  \left(\sum_{t=1}^T|\Delta u_1^n(t)-\Delta u_1^*(t)|^\theta+\sum_{t=1}^T|u_1^n(t)-u_1^*(t)|^\theta\right)^{1/\theta}\\
  &     &  +\left(\sum_{t=1}^T|\Delta  u_2^n(t)-\Delta  u_2^*(t)|^\theta+\sum_{t=1}^T|u_2^n(t)-u_2^*(t)|^\theta\right)^{1/\theta}\\
  &  =  &     \|u^n-u^*\|\to 0,\mbox{ as }n\to\infty,
 \end{eqnarray*}
 which implies that $\lim\limits_{n\to \infty}|\Delta u_i^n(t)-\Delta u_i^*(t)|=0$ and $\lim\limits_{n\to \infty}|u_i^n(t)-u_i^*(t)|=0, i=1,2$ for every $t\in \mathbb Z[1,T]$.
  Hence, it is easy to obtain that ($\mathcal{A}1$) implies that $I$ is continuous in $E$ and then sequentially weakly lower semicontinuous. Moreover, obviously, ($\mathcal{A}1$) and
  Lemma 2.1 imply that $I$  is a $C^1$ functional and
 \begin{eqnarray*}
 \langle I'(u),v\rangle&  =  &  \sum\limits_{t = 1}^{T}\left[\rho_1(t)(\phi_1(\Delta u_1(t)),\Delta v_1(t))+\rho_2(t)(\phi_2(\Delta u_2(t)),\Delta v_2(t))\right.\nonumber\\
             &     &  +\left.\rho_3(t)(\phi_3(u_1(t)), v_1(t))+\rho_4(t)(\phi_4(u_2(t)), v_2(t))\right],\quad \mbox{for } u,v\in E.
 \end{eqnarray*}
 It follows from ($\mathcal{A}3$) that
 \begin{eqnarray*}
          \langle I'(u)-I'(v),u-v\rangle
 &  =  &  \sum\limits_{t = 1}^{T}\left[\rho_1(t)(\phi_1(\Delta u_1(t))-\phi_1(\Delta v_1(t)),\Delta u_1(t)-\Delta v_1(t))\right.\nonumber\\
 &     &  +\rho_2(t)(\phi_2(\Delta u_2(t))-\phi_2(\Delta v_2(t)),\Delta u_2(t)-\Delta v_2(t))\nonumber\\
 &     &  +\rho_3(t)(\phi_3(u_1(t))-\phi_3(v_1(t)), u_1(t)-v_1(t))\nonumber\\
 &     & \left. +\rho_4(t)(\phi_4(u_1(t))-\phi_4(v_2(t)), u_2(t)-v_2(t))\right]\nonumber\\
 & \ge &  \sum\limits_{t = 1}^{T}\left[c_1\rho_1^- |\Delta u_1(t)-\Delta v_1(t)|^\theta+c_2\rho_2^-|\Delta u_2(t)-\Delta v_2(t)|^\theta\right.\nonumber\\
 &     & \quad\quad +\left.c_3 \rho_3^- |u_1(t)- v_1(t)|^\theta+c_4 \rho_4^-|u_2(t)- v_2(t)|^\theta\right]\nonumber\\
 & \ge & \min\{c_1 \rho_1^-,c_3\rho_3^-\}\|u_1-v_1\|_{E_T}^\theta+\min\{c_2\rho_2^-,c_4\rho_4^-\}\|u_2-v_2\|_{E_T}^\theta\nonumber\\
 & \ge & \frac{1}{2^{\theta-1}}\min\{c_1\rho_1^-,c_2\rho_2^-,c_3\rho_3^-,c_4\rho_4^-\}\left(\|u_1-v_1\|_{E_T}+\|u_2-v_2\|_{E_T}\right)^{\theta}\nonumber\\
 &  =  & \frac{1}{2^{\theta-1}}\min\{c_1\rho_1^-,c_2\rho_2^-,c_3\rho_3^-,c_4\rho_4^-\}\|u-v\|^{\theta}, \quad \mbox{for } u,v\in E.
 \end{eqnarray*}
 So $I'$ is uniformly monotone in $E$.
 By ($\mathcal{A}1$) and ($\mathcal{A}3$), we have
 \begin{eqnarray}
  \label{dd4}(\phi_i(x),x)\ge c_i|x|^\theta, \quad \mbox{for all } x\in \mathbb R^N, \quad  i=1,2,3,4.
 \end{eqnarray}
  Hence, (\ref{dd4}) implies that
  \begin{eqnarray}
  \label{dd3}
 &     &  \frac{ \langle I'(u),u\rangle}{\|u\|} \nonumber\\
 &  =  & \frac{1}{\|u\|}\sum\limits_{t = 1}^{T}\left[\rho_1(t)(\phi_1(\Delta u_1(t)),\Delta u_1(t))+\rho_2(t)(\phi_2(\Delta u_2(t)),\Delta u_2(t))\right.\nonumber\\
 &     &  +\left.\rho_3(t)(\phi_3(u_1(t)), u_1(t))+\rho_4(t)(\phi_4(u_2(t)), u_2(t))\right]\nonumber\\
 & \ge & \frac{1}{\|u\|}\left\{\sum\limits_{t = 1}^{T}\left[c_1\rho_1^-|\Delta u_1(t)|^\theta+c_2\rho_2^-|\Delta u_2(t)|^\theta+c_3\rho_3^-|u_1(t)|^\theta+c_4\rho_4^-|u_2(t)|^\theta\right]\right\}\nonumber\\
 & \ge &  \min\{c_1\rho_1^-,c_2\rho_2^-,c_3\rho_3^-,c_4\rho_4^-\}\frac{\|u_1\|_{E_T}^\theta+\|u_2\|_{E_T}^\theta}{\|u_1\|_{E_T}+\|u_2\|_{E_T}}\nonumber\\
 & \ge & \frac{1}{2^{\theta-1}}\min\{c_1\rho_1^-,c_2\rho_2^-,c_3\rho_3^-,c_4\rho_4^-\}\frac{(\|u_1\|_{E_T}+\|u_2\|_{E_T})^\theta}{\|u_1\|_{E_T}+\|u_2\|_{E_T}}\nonumber\\
 &  =  &  \frac{1}{2^{\theta-1}}\min\{c_1\rho_1^-,c_2\rho_2^-,c_3\rho_3^-,c_4\rho_4^-\}\|u\|^{\theta-1}
 \end{eqnarray}
 for all $u\in E$.  So $\lim\limits_{\|u\|\to\infty} \frac{ \langle I'(u),u\rangle}{\|u\|} =+\infty$, that is, $I'$ is coercive in $E$.
 Next, we show that $I'$ is also hemicontinuous in $E$. Assume that
 $s\to s^*$, $s,s^*\in [0,1]$. Note that
  \begin{eqnarray}
  \label{dd5}
   | \langle I'(u+sv),w\rangle-\langle I'(u+s^*v),w\rangle |
  \le   \|I'(u+sv)-I'(u+s^*v)\|\|w\|
 \end{eqnarray}
 for all $u,v,w\in E$. Then the continuity of $I'$  implies
 that $\langle I'(u+sv),w\rangle\to \langle I'(u+s^*v),w\rangle $ as
 $s\to s^*$ for all $u,v,w\in E$. Hence, $I'$ is hemicontinuous in $E$. Thus by Theorem
 26.A in \cite{Zeidler}, we know that $I'$ admits a continuous
 inverse in $E$.
 \par
 Obviously, ($\mathcal{A}2$) implies that $\Psi$, $\Phi$ and $\Gamma$ are
 $C^1$ functionals. Next, we show that $\Psi'$, $\Phi'$ and $\Gamma'$ are
 compact. Assume that $\{u^n\}\subset E$ is bounded. Then there
 exists a constant $D_1$ such that $\|u^n\|\le D_1$ and there exists a subsequence, still denoted by $\{u_n\}$, such that $u^n\rightharpoonup u^*$ for some
 $u^*\in E$. Furthermore, $u^n\to u^*$. By the continuity of $\Psi'$, $\Phi'$ and $\Gamma'$, it is clear that
  \begin{eqnarray*}
 \|\Psi'(u^n)-\Psi'(u^*)\|\to 0,\quad\|\Phi'(u^n)-\Phi'(u^*)\|\to 0,\quad\|\Gamma'(u^n)-\Gamma'(u^*)\|\to 0,\quad\mbox{ as } n\to\infty.
 \end{eqnarray*}
 Hence, $\Psi'$, $\Phi'$ and $\Gamma'$ are compact in $E$.
 It follows from ($\mathcal{A}5$) that
  \begin{eqnarray*}
  \Psi(u)+\lambda\Phi(u)=\sum_{t=1}^T\left[\lambda
  G(t,u_1(t),u_2(t))-F(t,u_1(t),u_2(t))\right]\ge TC_0(\lambda),
  \end{eqnarray*}
 which shows that $\Psi+\lambda\Phi$ is bounded below for all $\lambda>0$. Moreover,
 ($\mathcal{A}5$)  implies that for any positive constant $D_1$, there exists a positive constant $D_2(D_1)$, which depends on $D_1$, such that
 \begin{eqnarray}
 \label{fff5} F(t,x_1,x_2)\ge D_1(|x_1|^l+|x_2|^l)-D_2(D_1)
 \end{eqnarray}
 for all $(x_1,x_2)\in\mathbb R^N\times\mathbb R^N$ and $t\in\mathbb{Z}[1,T]$.
 Then (\ref{fff5}), ($\mathcal{A}4$), (\ref{d3})  and (\ref{d4})  imply that
  \begin{eqnarray*}
  &     &  \lim_{\|u\|\to\infty}\frac{\Psi(u)}{I(u)}\\
  &  =  &  \lim_{\|u\|\to\infty}\frac{- \sum\limits_{t = 1}^{T}  F(t,u_{1}(t),u_{2}(t))}{\sum\limits_{t = 1}^{T}\left[\rho_1(t)\Phi_1(\Delta u_{1}(t))+\rho_2(t)\Phi_2(\Delta  u_{2}(t))+\rho_3(t)\Phi_3(u_1(t))+\rho_4(t)\Phi_4(u_2(t))\right]}\\
  & \le &  \lim_{\|u\|\to\infty}\frac{- D_1\sum\limits_{t = 1}^{T} (|u_1(t)|^l+|u_2(t)|^l)+D_2(D_1)T}{\sum\limits_{t = 1}^{T}\left[\rho_1(t)\Phi_1(\Delta u_{1}(t))+\rho_2(t)\Phi_2(\Delta  u_{2}(t))+\rho_3(t)\Phi_3(u_1(t))+\rho_4(t)\Phi_4(u_2(t))\right]}\\
  & \le &  \lim_{\|u\|\to\infty}\frac{- D_1\sum\limits_{t = 1}^{T} (|u_1(t)|^l+|u_2(t)|^l)}{\sum\limits_{t = 1}^{T}\left[d_1\rho_1^+|\Delta u_{1}(t)|^l+d_2\rho_2^+|\Delta  u_{2}(t)|^l+d_3\rho_3^+|u_1(t)|^l+d_4\rho_4^+|u_2(t)|^l+\sum\limits_{i=1}^4m_i\rho_i^+\right]}\\
  &     &  +\lim_{\|u\|\to\infty}\frac{D_2(D_1)T}{\min\{\rho_1^-,\rho_2^-,\rho_2^-,\rho_4^-\}\sum\limits_{t = 1}^{T}\left[\Phi_1(\Delta u_{1}(t))+\Phi_2(\Delta  u_{2}(t))+\Phi_3(u_1(t))+\Phi_4(u_2(t))\right]}\\
  & \le &  \lim_{\|u\|\to\infty}\frac{-D_1C_3^l(\|u_1\|_{E_T}^l+\|u_2\|_{E_T}^l)}{\max\{d_1\rho_1^+,d_2\rho_2^+,d_3\rho_3^+,d_4\rho_4^+\}(\|u_1\|_{[E_T]}^l+\|u_2\|_{[E_T]}^l)+\sum\limits_{t = 1}^{T}\sum\limits_{i=1}^4m_i\rho_i^+}\\
  & \le &  \lim_{\|u\|\to\infty}\frac{-D_1C_3^l\frac{1}{2^{l-1}}(\|u_1\|_{E_T}+\|u_2\|_{E_T})^l}{\max\{d_1\rho_1^+,d_2\rho_2^+,d_3\rho_3^+,d_4\rho_4^+\}(\|u_1\|_{[E_T]}^l+\|u_2\|_{[E_T]}^l)+\sum\limits_{t = 1}^{T}\sum\limits_{i=1}^4m_i\rho_i^+}\\
  & \le &  \lim_{\|u\|\to\infty}\frac{-D_1C_3^l\frac{1}{2^{l-1}}(\|u_1\|_{E_T}+\|u_2\|_{E_T})^l}{\max\{d_1\rho_1^+,d_2\rho_2^+,d_3\rho_3^+,d_4\rho_4^+\}(C_{6}^l\|u_1\|_{E_T}^l+C_{6}^l\|u_2\|_{E_T}^l)+\sum\limits_{t = 1}^{T}\sum\limits_{i=1}^4m_i\rho_i^+}\\
  & \le &  \lim_{\|u\|\to\infty}\frac{-D_1C_3^l\frac{1}{2^{l-1}}(\|u_1\|_{E_T}+\|u_2\|_{E_T})^l}{\max\{d_1\rho_1^+,d_2\rho_2^+,d_3\rho_3^+,d_4\rho_4^+\}C_{6}^l(\|u_1\|_{E_T}+\|u_2\|_{E_T})^l+\sum\limits_{t = 1}^{T}\sum\limits_{i=1}^4m_i\rho_i^+}\\
  &  =  &  \frac{1}{2^{l-1}}\frac{-D_1C_3^l}{\max\{d_1\rho_1^+,d_2\rho_2^+,d_3\rho_3^+,d_4\rho_4^+\}C_{6}^l}.
  \end{eqnarray*}
 By the arbitrary of $D_1$, we obtain that
  \begin{eqnarray*}
     \lim_{\|u\|\to\infty}\frac{\Psi(u)}{I(u)}=-\infty.
  \end{eqnarray*}
  By ($\mathcal{A}1$) and ($\mathcal{A}4$), we know that $\Phi_i$
  reaches its unique minimum at $0$, $i=1,2,3,4$ (see \cite{Mawhin2012}) and so $I$ has unique global minima
  $0$. Then $M=\{0\}$. By ($\mathcal{A}6$), we have $\sup_M\Phi=0$.  Hence, by Theorem 2.1, the conclusion of Theorem 1.1 holds. $\Box$

 \vskip2mm
 \noindent
 {\bf Proof of Corollary 1.1}\ \  It follows from ($\mathcal{A}5$)$'$ that there exist $D_3>0$ and  $D_4>0$ such
 that for every $t\in \mathbb Z[1,T]$,
 \begin{eqnarray*}
 \label{fff2} F(t,x_1,x_2)\le D_3|x_1|^s+D_3|x_2|^s+D_4
 \end{eqnarray*}
 and for any  $D_5>D_3$, there are a constant
 $D_6(D_5)$, which depends on $D_5$, such that
 \begin{eqnarray*}
 \label{fff4}  G(t,x_1,x_2)\ge D_5|x_1|^s+D_5|x_2|^s+D_6(D_5).
 \end{eqnarray*}
 Obviously, for every $\lambda>0$, we can find a sufficiently large $D_5(\lambda)$ such that $\lambda
 D_5(\lambda)>D_3$. Hence, we have
 \begin{eqnarray*}
 \label{fff4}  \lambda G(t,x_1,x_2)\ge D_3|x_1|^s+D_3|x_2|^s+\lambda D_6(D_5(\lambda))\ge F(t,x_1,x_2)-D_4+\lambda D_6(D_5(\lambda)).
 \end{eqnarray*}
 So ($\mathcal{A}5$)$'$ implies ($\mathcal{A}5$). $\Box$

  \section{Proof of Theorem 1.2}
 \setcounter{equation}{0}
 \par
 \quad\quad When the condition ($\mathcal{\gamma}$) holds, on $E_T$, we define
  \begin{eqnarray*}
     \|u\|_{(E_{T,q})}=\left(\sum_{t=1}^T\gamma_1(t)|\Delta   u_1(t)|^q+\sum_{t=1}^T\gamma_3(t)|u_1(t)|^q\right)^{1/q}
  \end{eqnarray*}
  and
  \begin{eqnarray*}
        \|u\|_{(E_{T,p})}=\left(\sum_{t=1}^T\gamma_2(t)|\Delta   u_2(t)|^p+\sum_{t=1}^T\gamma_4(t)|u_2(t)|^p\right)^{1/p}.
  \end{eqnarray*}
 For $u=(u_1,u_2)\in E$, define
 $$
 \|u\|_{(\infty)}=\|u_1\|_{\infty}+\|u_2\|_{\infty}.
 $$
 Moreover, it is clear that $E$ is homeomorphic to $\mathbb R^{2NT}$. Then there is a basis of $E$ denoted by $\{e_1, e_2,..., e_{2NT}\}$. For every $u\in E$, there exists a unique point $(\lambda_1,\lambda_2,...,\lambda_{2NT})\in \mathbb{R}^{2NT}$ such that
  \begin{eqnarray*}
  u=\sum_{i=1}^{2NT}\lambda_ie_i
  \end{eqnarray*}
  and define
  $$\|u\|_{(2)}=\left(\sum_{i=1}^{2NT}\lambda_i^2\right)^{\frac{1}{2}}.$$
  Set
  $$E_{\delta}=\{u\in E: \|u\|_{(2)}=\delta\}.$$
 Since both $E$ and $E_T$ are finite-dimensional spaces, then $\|\cdot\|_{(\infty)}$ is equivalent to $\|\cdot\|_{(2)}$ on $E$, and both $\|\cdot\|_{(E_T,q)}$  and  $\|\cdot\|_{(E_T,p)}$  are equivalent to
  $\|\cdot\|_{\infty}$ on $E_T$. Hence, there exist positive constants $R_i$ $(i=1,2,\cdots,6)$ such that
 \begin{eqnarray}
 &   & \label{dvv1}  R_1\|\cdot\|_{(2)}\le \|\cdot\|_{(\infty)}\le R_2\|\cdot\|_{(2)},\\
 &   & \label{dv1}  R_3\|\cdot\|_{\infty}\le \|\cdot\|_{(E_{T,q})}\le R_4\|\cdot\|_{\infty},\\
 &   & \label{dv2}  R_5\|\cdot\|_{\infty}\le \|\cdot\|_{(E_{T,p})}\le R_6\|\cdot\|_{\infty}.
 \end{eqnarray}
 \par
 In Lemma 2.1,  let
 $$
 L(t,x_{1},x_{2},y_{1},y_{2})=\gamma_1(t)\Phi_1(y_{1})+\gamma_2(t)\Phi_2(y_{2})+\gamma_3(t)\Phi_3(x_1)+\gamma_4(t)\Phi_4(x_2)- F(t,x_{1},x_{2}),
 $$
 where $ F $: $\mathbb{Z}[1,T]\times\mathbb{R}^{N}\times\mathbb{R}^{N}\longrightarrow\mathbb{R}$ is continuously differentiable in
 $(x_{1},x_{2})\in \mathbb{R}^{N}\times\mathbb{R}^{N}$ for all $t\in\mathbb{Z}[1,T]$.
 Then
 \begin{eqnarray}
 \label{bv1}         \varphi(u)
 &  =  &  \sum\limits_{t = 1}^{T}\gamma_1(t)\Phi_1(\Delta u_1(t))+\sum\limits_{t = 1}^{T}\gamma_2(t)\Phi_2(\Delta u_2(t))\nonumber\\
 &     &  +\sum\limits_{t = 1}^{T}\gamma_3(t)\Phi_3(u_1(t))+\sum\limits_{t = 1}^{T}\gamma_4(t)\Phi_4(u_2(t))-\sum\limits_{t=1}^{T}F(t,u_{1}(t),u_{2}(t)).
 \end{eqnarray}
 and for $ \forall\ u,v\in E$, we have
 \begin{eqnarray}
  \label{ffv2}&     &    \langle\varphi'(u),v\rangle\nonumber\\
            &  =  &  \langle\varphi'(u_{1},u_{2}),(v_{1},v_{2})\rangle\nonumber\\
            &  =  &  \sum\limits_{t = 1}^{T}\left[\gamma_1(t)(\phi_1(\Delta u_1(t)),\Delta v_1(t))+\gamma_2(t)(\phi_2(\Delta u_2(t)),\Delta v_2(t))\right.\nonumber\\
            &     &  \left.+\gamma_3(t)(\phi_3( u_1(t)), v_1(t))+ \gamma_4(t)(\phi_4( u_2(t)), v_2(t))\right.\nonumber\\
            &     &  -\sum\limits_{t=1}^{T}(\nabla_{u_1} F(t,u_{1}(t),u_{2}(t)),v_1(t))-\sum\limits_{t=1}^{T}(\nabla_{u_2}F(t,u_{1}(t),u_{2}(t)),v_2(t)).
  \end{eqnarray}

 \vskip2mm
 \par
 Similar to the argument of Lemma 2.2, it is easy to obtain the
 following Lemma:
 \vskip2mm
 \noindent
 {\bf Lemma 4.1} {\it If $u\in E$ is a solution of Euler equation
 $\varphi'(u)=0$, then  $u$ is a solution of system (\ref{eqq1}).}

 \vskip2mm
 \par
Denote with $\theta$ the zero element of $X$ and with $\Sigma$  the
 family of sets $A\subset X\backslash \{\theta\}$ such that $A$ is
closed in $X$ and symmetric with respect to $\theta,$ i.e. $u\in A$
implies $-u\in A.$
 \vskip2mm
 \noindent
 {\bf Theorem 4.1} (see \cite{Ra}, Theorem 9.1) \ \ {\it Let $X$ be a real
Banach space and  $\varphi$ be an even function belonging to
$C^1(X,\mathbb{R})$ with $\varphi(\theta)=0$, bounded from below and
satisfying (PS) condition. Suppose that there is a set $K\in \Sigma$
such that $K$ is  homeomorphic to $S^{j-1}$($j-1$ dimension unit
sphere) by an odd map and $\sup_K \varphi<0.$ Then, $\varphi$ has at
least $j$ distinct pairs of nonzero critical points.}

 \vskip2mm
 \noindent
 {\bf Proof of Theorem 1.2}\ \ It follows from $(\phi)$, (\ref{bv1}), (\ref{ff1}), (\ref{dv1}), (\ref{dv2}) and ($\mathcal{F}$0) that
  \begin{eqnarray}
 \label{cccv1}         \varphi(u)
 &  =  &  \sum\limits_{t = 1}^{T}\gamma_1(t)\Phi_1(\Delta u_{1}(t))+\sum\limits_{t = 1}^{T}\gamma_2(t)\Phi_2(\Delta u_{2}(t))\nonumber\\
 &     &  +\sum\limits_{t = 1}^{T}\gamma_3(t)\Phi_3(u_{1}(t))+\sum\limits_{t = 1}^{T}\gamma_4(t)\Phi_4(u_{2}(t))-\sum\limits_{t=1}^{T}F(t,u_{1}(t),u_{2}(t))\nonumber\\
 & \ge &  a_1\sum\limits_{t = 1}^{T}\gamma_1(t)|\Delta u_{1}(t)|^q+a_2\sum\limits_{t = 1}^{T}\gamma_2(t)|\Delta u_{2}(t)|^p\nonumber\\
 &     &  +a_3\sum\limits_{t = 1}^{T}\gamma_3(t)|u_{1}(t)|^q+a_4\sum\limits_{t = 1}^{T}\gamma_4(t)|u_{2}(t)|^p-\sum\limits_{t=1}^{T}F(t,u_{1}(t),u_{2}(t))\nonumber\\
 & \ge &   \min\{a_1,a_3\}\|u_1\|_{(E_{T,q})}^q+ \min\{a_2,a_4\}\|u_2\|_{(E_{T,p})}^p\nonumber\\
 &     &  -\sum\limits_{t=1}^{T}\left[h_1(t)|u_1(t)|^{\alpha_1}+h_2(t)|u_2(t)|^{\alpha_2}+l(t)\right]\nonumber\\
 & \ge &   \min\{a_1,a_3\}R_{3}^{q}\|u_1\|_{\infty}^q+ \min\{a_2,a_4\}R_{5}^{p}\|u_2\|_{\infty}^p\nonumber\\
 &     &  - \|u_1\|_{\infty}^{\alpha_1}\sum\limits_{t=1}^{T}h_1(t)-\|u_2\|_{\infty}^{\alpha_2}\sum\limits_{t=1}^{T}h_2(t) -\sum\limits_{t=1}^{T}l(t)
 \end{eqnarray}
 for all $u\in E$. Since $\alpha_1\in [0,q)$ and $\alpha_2\in
  [0,p)$, it is easy to see that
  \begin{eqnarray}
   \label{cv4}\varphi(u)\to+\infty, \mbox{ as }\|u\|_{(\infty)}=\|u_1\|_{\infty}+\|u_2\|_{\infty}\to\infty,
  \end{eqnarray}
  which implies that $\varphi$ is bounded from below and any (PS) sequence $\{u_n\}$ is
  bounded. Hence $\varphi$ satisfies (PS) condition.
  Obviously, ($\mathcal{F}$1) and ($\mathcal{F}$2)
 imply that $\varphi(0)=0$ and $\varphi$ is even. Next, we prove
 that there exists a set $K\subset E$ such that $K$ is homeomorphic
 to $S^{2NT-1}$ by an odd map, and $\sup_K \varphi<0$. Note that $\delta<1$.
 For all $u=(u_1,u_2)\in E_\delta$ and $r>0$, by (\ref{dvv1}) we have
 \begin{eqnarray}\label{cv2}
 &     & M_1r^{\beta_1}\|u_1\|_\infty^{\beta_1}+M_2r^{\beta_2}\|u_2\|_\infty^{\beta_2}\nonumber\\
 &  =  & M_1r^{\beta_1}R_2^{\beta_1}\left\|\frac{u_1}{R_2}\right\|_\infty^{\beta_1}+M_2r^{\beta_2}R_2^{\beta_2}\left\|\frac{u_2}{R_2}\right\|_\infty^{\beta_2}\nonumber\\
 & \ge & \min\{M_1r^{\beta_1}R_2^{\beta_1},M_2r^{\beta_2}R_2^{\beta_2}\}
         \left(\left\|\frac{u_1}{R_2}\right\|_\infty^{\max\{\beta_1,\beta_2\}}+\left\|\frac{u_2}{R_2}\right\|_\infty^{\max\{\beta_1,\beta_2\}}\right)\nonumber\\
 & \ge & 2^{1-\max\{\beta_1,\beta_2\}}\min\{M_1r^{\beta_1}R_2^{\beta_1},M_2r^{\beta_2}R_2^{\beta_2}\}
         \left(\left\|\frac{u_1}{R_2}\right\|_\infty+\left\|\frac{u_2}{R_2}\right\|_\infty\right)^{\max\{\beta_1,\beta_2\}}\nonumber\\
 &  =  & 2^{1-\max\{\beta_1,\beta_2\}}\min\{M_1r^{\beta_1}R_2^{\beta_1},M_2r^{\beta_2}R_2^{\beta_2}\}\left(\frac{1}{R_2}\right)^{\max\{\beta_1,\beta_2\}}
         \|u\|_{(\infty)}^{\max\{\beta_1,\beta_2\}}\nonumber\\
 & \ge & 2^{1-\max\{\beta_1,\beta_2\}}\min\{M_1r^{\beta_1}R_2^{\beta_1},M_2r^{\beta_2}R_2^{\beta_2}\}\left(\frac{1}{R_2}\right)^{\max\{\beta_1,\beta_2\}}
         R_1^{\max\{\beta_1,\beta_2\}}\|u\|_{(2)}^{\max\{\beta_1,\beta_2\}}\nonumber\\
 &  =  & 2\min\{M_1r^{\beta_1}R_2^{\beta_1},M_2r^{\beta_2}R_2^{\beta_2}\}\left(\frac{R_1\delta}{2R_2}\right)^{\max\{\beta_1,\beta_2\}}.
  \end{eqnarray}
  Then for all $u=(u_1,u_2)\in E_\delta$ and $0<r<\frac{1}{R_2}$, by $(\phi)$, ($\mathcal{F}$3), (\ref{ff1}), (\ref{dvv1})-(\ref{dv2}) and (\ref{cv2}) we have
   \begin{eqnarray}
 \label{cv1}         \varphi(ru)
 &  =  &  \sum\limits_{t = 1}^{T}\gamma_1(t)\Phi_1(r\Delta u_{1}(t))+\sum\limits_{t = 1}^{T}\gamma_2(t)\Phi_2(r\Delta u_{2}(t))\nonumber\\
 &     &  +\sum\limits_{t = 1}^{T}\gamma_3(t)\Phi_3(ru_{1}(t))+\sum\limits_{t = 1}^{T}\gamma_4(t)\Phi_4(ru_{2}(t))-\sum\limits_{t=1}^{T}F(t,ru_{1}(t),ru_{2}(t))\nonumber\\
 & \le &  b_1\sum\limits_{t = 1}^{T}\gamma_1(t)|r\Delta u_{1}(t)|^q+b_2\sum\limits_{t = 1}^{T}\gamma_2(t)|r\Delta u_{2}(t)|^p\nonumber\\
 &     &  +b_3\sum\limits_{t = 1}^{T}\gamma_3(t)|ru_{1}(t)|^q+b_4\sum\limits_{t = 1}^{T}\gamma_4(t)|ru_{2}(t)|^p-\sum\limits_{t=1}^{T}F(t,ru_{1}(t),ru_{2}(t))\nonumber\\
 & \le & \max\{b_1,b_3\}r^q\|u_1\|_{(E_T,q)}^q+\max\{b_2,b_4\}r^p\|u_2\|_{(E_T,p)}^p\nonumber\\
 &     &    -M_1r^{\beta_1}\sum\limits_{t=1}^{T}|u_{1}(t)|^{\beta_1}-M_2r^{\beta_2}\sum\limits_{t=1}^{T}|u_{2}(t)|^{\beta_2}\nonumber\\
 & \le & \max\{b_1,b_3\}r^qR_4^q\|u_1\|_{\infty}^q+\max\{b_2,b_4\}r^pR_6^p\|u_2\|_{\infty}^p\nonumber\\
 &     &    -M_1r^{\beta_1}\|u_{1}\|_{\infty}^{\beta_1}-M_2r^{\beta_2}\|u_{2}\|_{\infty}^{\beta_2}\nonumber\\
 & \le & \max\{b_1,b_3\}r^qR_4^qR_2^q\|u\|_{(2)}^q+\max\{b_2,b_4\}r^pR_6^pR_2^p\|u\|_{(2)}^p\nonumber\\
 &     &    -2\min\{M_1r^{\beta_1}R_2^{\beta_1},M_2r^{\beta_2}R_2^{\beta_2}\}\left(\frac{R_1\delta}{2R_2}\right)^{\max\{\beta_1,\beta_2\}}\nonumber\\
 &  =  & \max\{b_1,b_3\}r^q(R_4R_2\delta)^q+\max\{b_2,b_4\}r^p(R_6R_2\delta)^p\nonumber\\
 &     &    -2\min\{M_1r^{\beta_1}R_2^{\beta_1},M_2r^{\beta_2}R_2^{\beta_2}\}\left(\frac{R_1\delta}{2R_2}\right)^{\max\{\beta_1,\beta_2\}}.
 \end{eqnarray}
 Since $\beta_i\in (1,\min\{q,p\})$, $i=1,2$. Then (\ref{cv1}) implies that there exist sufficiently small $r_0\in (0,1)$ and $\epsilon>0$ such
 that there exists sufficiently small $r_0\in (0,1)$ and $\epsilon>0$ such that $\varphi(r_0 u)<-\epsilon$ for all $u\in E_\delta$. Set
 $$E_{\delta}^{r_0}=\{r_0u: u\in E_{\delta}\}\quad \mbox{ and }\quad S^{2NT-1}=\left\{(\lambda_1,\lambda_2,\cdots,\lambda_{2NT})\in \mathbb R^{2NT}: \sum_{i=1}^{2NT}  \lambda_i^2=1\right\}.$$
 Then $E_\delta^{r_0}\in \Sigma$ and
 \begin{eqnarray}\label{cv3}
 \psi(u)<-\epsilon, \quad\forall\ u\in E_\delta^{r_0}.
 \end{eqnarray}
 Define the map $\psi:E_\delta^{r_0}\to S^{2NT-1}$ by
 $$
 \psi(u)=\psi\left(\sum_{i=1}^{2NT}\lambda_ie_i\right)=\frac{1}{r_0\delta}(\lambda_1,\lambda_2,\cdots,\lambda_{2NT}).
 $$
 Then it is easy to see that $\psi$ is an odd and homeomorphic map.
 Moreover, (\ref{cv3}) implies that $\sup_{E_\delta^{r_0}}\varphi\le-\epsilon<0$.
 Therefore, by Theorem 4.1, we obtain that system (\ref{eqq1}) has at least
 $2NT$ distinct pairs of solutions in $E$. $\Box$

 \section{Examples}
 \setcounter{equation}{0}
 \vskip1mm
 \noindent
 {\bf Example 5.1.} We present an example to which Theorem 1.1 applies and make an estimate for the parameters in our result. Let $T=2$ and $N$ be fixed integer. Assume that $\phi_1(y)=y+|y|^{\frac{1}{3}}y$,
 $\phi_2(y)=y+ |y|^{\frac{1}{2}}y$,
 $\phi_3(y)=\phi_4(y)=2y$,
 $\rho_i$ are $2$-periodic and satisfy
 $\rho_i> 0$ for all $t\in \mathbb Z[1,2]$ ($i=1,2,3,4$).
 Then  $\Phi_1(y)=\frac{|y|^{2}}{2}+ \frac{|y|^{\frac{7}{3}}}{\frac{7}{3}}$,  $\Phi_2(y)=\frac{|y|^{2}}{2}+ \frac{|y|^{\frac{5}{2}}}{\frac{5}{2}}$, $\Phi_3(y)=\Phi_4(y)=|y|^{2}$. Let
  \begin{eqnarray}
  \label{ddd1}         F(t,x_1,x_2)
 &  =  &  |x_1|^3+|x_2|^3\nonumber\\
 G (t,x_1,x_2)& =  & |x_1|^4+|x_2|^4\nonumber\\
 H (t,x_1,x_2)& =  & (\cos^2\frac{\pi t}{2}+2)\sin(|x_1|^2+|x_2|^2+2)\nonumber\\
 W(t,x_1,x_2) & =  & F(t,x_1,x_2)-\lambda G (t,x_1,x_2)+\nu H (t,x_1,x_2).\nonumber
 \end{eqnarray}
 Then
 \begin{eqnarray*}
  &   &  I(u)   =    \sum\limits_{t = 1}^{2}\left[\rho_1(t)\left(\frac{|\Delta u_{1}(t)|^{2}}{2}+ \frac{|\Delta u_{1}(t)|^{\frac{7}{3}}}{\frac{7}{3}}\right)+\rho_2(t)\left(\frac{|\Delta u_{2}(t)|^{2}}{2}+ \frac{|\Delta u_{2}(t)|^{\frac{5}{2}}}{\frac{5}{2}}\right)\right]\\
  &&\qquad\quad+\sum\limits_{t = 1}^{2}\left[\rho_3(t)|u_1(t)|^2+\rho_4(t)|u_2(t)|^2\right],\\
  &   &  \Psi(u)   =   -  \sum\limits_{t = 1}^{2}  \left(|u_1(t)|^3+|u_2(t)|^3\right),\quad \Phi(u)   =     \sum\limits_{t = 1}^{2}  \left(|u_1(t)|^4+|u_2(t)|^4\right),\\
   &   &   \Gamma(u)   =   -\sum\limits_{t = 1}^{T}  (\cos^2\frac{\pi t}{2}+2)\sin(|u_1(t)|^2+|u_2(t)|^2+2), \quad\quad u\in   E,
 \end{eqnarray*}
 where the definition of $E$ and its norm are in section 2.
 Take $\theta=2$ and $l=\frac{5}{2}$. With a similar discussion as in Remark 1.1, we can prove that  all conditions of Theorem 1.1 hold. Since, by a simply computation we have
 \begin{eqnarray}\label{5.1}
 \gamma=\inf_E(\Psi(u)+\Phi(u))=\inf_E\sum\limits_{t = 1}^{2}  \left(|u_1(t)|^4+|u_2(t)|^4-|u_1(t)|^3-|u_2(t)|^3\right)=-\frac{27}{64},
 \end{eqnarray}
 which is obtained when $|u_1(1)|=|u_1(2)|=|u_2(1)|=|u_2(2)|=\frac{3}{4}.$ Moreover, for $r>0$, we have
 \begin{eqnarray*}
 \Phi^{-1}(r)          &=&\{u\in E: |u_1(1)|^4+|u_1(2)|^4+|u_2(1)|^4+|u_2(2)|^4=r\},\\
 \Phi^{-1}(]-\infty,r])&=&\{u\in E: |u_1(1)|^4+|u_1(2)|^4+|u_2(1)|^4+|u_2(2)|^4\le r\},\\
 \Phi^{-1}(]r,+\infty[)&=&\{u\in E: |u_1(1)|^4+|u_1(2)|^4+|u_2(1)|^4+|u_2(2)|^4>r\}.
 \end{eqnarray*}
 Then
 \begin{eqnarray}\label{5.2}
 \eta_r & = & \inf_{u\in\Phi^{-1}(r)}I(u)\nonumber\\
        & = & \inf_{u\in\Phi^{-1}(r)}\left\{ \sum\limits_{t = 1}^{2}\left[\rho_1(t)\left(\frac{|\Delta u_{1}(t)|^{2}}{2}+ \frac{|\Delta u_{1}(t)|^{\frac{7}{3}}}{\frac{7}{3}}\right)+\rho_2(t)\left(\frac{|\Delta u_{2}(t)|^{2}}{2}+ \frac{|\Delta u_{2}(t)|^{\frac{5}{2}}}{\frac{5}{2}}\right)\right]\right.\nonumber\\
               &&\left.\qquad\qquad+\sum\limits_{t = 1}^{2}\left[\rho_3(t)|u_1(t)|^2+\rho_4(t)|u_2(t)|^2\right]\right\}\nonumber\\
        & \geq & \inf_{u\in\Phi^{-1}(r)}\sum\limits_{t = 1}^{2}\left[\rho_3(t)|u_1(t)|^2+\rho_4(t)|u_2(t)|^2\right]\nonumber\\
        &  =   & \inf_{u\in\Phi^{-1}(r)}\left(\rho_3(1)|u_1(1)|^2+\rho_3(2)|u_1(2)|^2+\rho_4(1)|u_2(1)|^2+\rho_4(2)|u_2(2)|^2\right)\nonumber\\
        &  =   & \min\{\rho_3(1),\rho_3(2),\rho_4(1),\rho_4(2)\}\sqrt{r},
 \end{eqnarray}
 which can be obtained by using the Lagrange multiplier method. By (\ref{5.1}), (\ref{5.2}) and the fact that $\Phi(0)=I(0)=0$, we have
 \begin{eqnarray}\label{5.3}
 \mu^*(I,\Psi,\Phi,r)&  =  & \inf\left\{\frac{\Psi(u)-\gamma+r}{\eta_r-I(u)}:u\in E,\Phi(u)<r, I(u)<\eta_r\right\}\nonumber\\
                     & \le & \frac{\Psi(0)-\gamma+r}{\eta_r-I(0)}=\frac{\frac{27}{64}+r}{\eta_r}
                             \le\frac{\frac{27}{64}+r}{\min\{\rho_3(1),\rho_3(2),\rho_4(1),\rho_4(2)\}\sqrt{r}}.
 \end{eqnarray}
 When $\mu>\frac{r^{\frac{1}{4}}}{\min\{\rho_3(1),\rho_3(2),\rho_4(1),\rho_4(2)\}}$, we have
 \begin{eqnarray}\label{5.4}
 & & \inf\limits_{\Phi^{-1}(]-\infty,r])}(\mu I+\Psi)\nonumber\\
 &=& \inf\limits_{\Phi^{-1}(]-\infty,r])}\mu\left\{ \sum\limits_{t = 1}^{2}\left[\rho_1(t)\left(\frac{|\Delta u_{1}(t)|^{2}}{2}+ \frac{|\Delta u_{1}(t)|^{\frac{7}{3}}}{\frac{7}{3}}\right)+\rho_2(t)\left(\frac{|\Delta u_{2}(t)|^{2}}{2}+ \frac{|\Delta u_{2}(t)|^{\frac{5}{2}}}{\frac{5}{2}}\right)\right]\right.\nonumber\\
               &&\left.\qquad\qquad+\sum\limits_{t = 1}^{2}\left[\rho_3(t)|u_1(t)|^2+\rho_4(t)|u_2(t)|^2\right]\right\}
               -  \sum\limits_{t = 1}^{2}  \left(|u_1(t)|^3+|u_2(t)|^3\right)\nonumber\\
 &=& 0,
 \end{eqnarray}
 which is obtained when $u_1(1)=u_1(2)=u_2(1)=u_2(2)=0.$ When $\mu>\frac{r^{\frac{1}{4}}}{\min\{\rho_3(1),\rho_3(2),\rho_4(1),\rho_4(2)\}}$, we choose $u_0: u_1(1)=u_1(2)=u_2(1)=u_2(2)=\mu[\rho_3(1)+\rho_3(2)+\rho_4(1)+\rho_4(2)]$, then
 \begin{eqnarray*}
 &   & |u_1(1)|^4+|u_1(2)|^4+|u_2(1)|^4+|u_2(2)|^4\\
 & = & 4\mu^4[\rho_3(1)+\rho_3(2)+\rho_4(1)+\rho_4(2)]^4\\
 & > & \frac{4r[\rho_3(1)+\rho_3(2)+\rho_4(1)+\rho_4(2)]^4}{(\min\{\rho_3(1),\rho_3(2),\rho_4(1),\rho_4(2)\})^4}\\
 & > & r
 \end{eqnarray*}
 which implies that $u_0\in\Phi^{-1}(]r,+\infty[).$
  Therefore, when $\mu>\frac{\max\{\frac{\frac{27}{64}+r}{\sqrt{r}}, r^{\frac{1}{4}}\}}{\min\{\rho_3(1),\rho_3(2),\rho_4(1),\rho_4(2)\}}$, by (\ref{5.4}) and the fact that $u_0\in\Phi^{-1}(]r,+\infty[)$ we have
 \begin{eqnarray}\label{5.5}
 \beta(\mu I+\Psi,\Phi,r)& = &\sup_{u\in\Phi^{-1}(]r,+\infty[)}\frac{\mu  I(u)+\Psi(u)-\inf\limits_{\Phi^{-1}(]-\infty,r])}(\mu I+\Psi)}{r-\Phi(u)}\nonumber\\
                         & = & \sup_{u\in\Phi^{-1}(]r,+\infty[)}\frac{\mu  I(u)+\Psi(u)}{r-\Phi(u)}\nonumber\\
                         &\ge&\frac{\mu  I(u_0)+\Psi(u_0)}{r-\Phi(u_0)}\nonumber\\
                         & = & \frac{3\mu^3[\rho_3(1)+\rho_3(2)+\rho_4(1)+\rho_4(2)]^3}{4\mu^4[\rho_3(1)+\rho_3(2)+\rho_4(1)+\rho_4(2)]^4-r}.
 \end{eqnarray}
 Hence, by Theorem 1.1, we obtain that for each $r>0$, for each $\mu>\frac{\max\{\frac{\frac{27}{64}+r}{\sqrt{r}}, r^{\frac{1}{4}}\}}{\min\{\rho_3(1),\rho_3(2),\rho_4(1),\rho_4(2)\}}$, and
 for each compact interval $[a,b]\subset]0,\frac{3\mu^3[\rho_3(1)+\rho_3(2)+\rho_4(1)+\rho_4(2)]^3}{4\mu^4[\rho_3(1)+\rho_3(2)+\rho_4(1)+\rho_4(2)]^4-r}[$, there exists a number $\rho>0$ with the following
 property: for every $\lambda\in [a,b]$, there exists
 $\delta>0$ such that, for each $\nu\in [0,\delta]$, system
  (\ref{eq1}) has at least three  $2$-periodic solutions in $E$ whose norms are less than $\rho$.

 \vskip1mm
 \noindent
 {\bf Example 5.2.} We present this example when Theorem
 1.2 applies.
  Let $N=6$ and $T=4$. Assume that $\phi_1(x)=\phi_3(x)=|x|^3x$, $\phi_2(x)=\phi_4(x)=|x|x$. Consider the following
  nonlinear difference
 ($\phi_1,\phi_2$)-Laplacian  system:
  \begin{eqnarray}\label{eqq3}
 \left\{
  \begin{array}{lll}
 \Delta\left(\left[\sin^2\frac{\pi}{4}(t-1)+1\right] \phi_{1}\big(\Delta u_{1}(t-1)\big)\right)-\left[|\sin\frac{\pi}{4}t|+1\right]\phi_{3}(|u_1(t)|)\\
 \quad\quad\quad\quad\quad\quad\quad\quad\quad\quad\quad\quad\quad\quad\quad\quad\quad\quad+\nabla_{u_{1}} F\big(t,u_{1}(t),u_{2}(t)\big)=0 \\
 \Delta\left(\left[\cos^2\frac{\pi}{4}(t-1)+1\right] \phi_{2}\big(\Delta u_{2}(t-1)\big)\right)-\left[|\cos\frac{\pi}{4}t|+1\right]\phi_{4}(|u_2(t)|)\\
 \quad\quad\quad\quad\quad\quad\quad\quad\quad\quad\quad\quad\quad\quad\quad\quad\quad\quad+\nabla_{u_{2}} F\big(t,u_{1}(t),u_{2}(t)\big)=0.
    \end{array}
 \right.
 \end{eqnarray}
 Then  $\gamma_1(t)=\sin^2\frac{\pi}{4}t+1$,
 $\gamma_2(t)=\cos^2\frac{\pi}{4}t+1$,
 $\gamma_3(t)=|\sin\frac{\pi}{4}t|+1$, $\gamma_4(t)=|\cos\frac{\pi}{4}t|+1 $.
 Obviously, the conditions ($\mathcal{\gamma}$) and $(\mathcal{\phi})$ hold and
 $\gamma_i$, $i=1,2,3,4$ are $T$-periodic ($T=4$).
 \par
 If we assume that
 $$
  F(t,x_1,x_2)=\left(\left|\sin\frac{\pi}{4}t\right|+1\right)|x_1|^{\frac{3}{2}}+\left(\cos^2\frac{\pi}{4}t+1\right)|x_2|^{2},
 $$
 then, obviously, ($\mathcal{F}$0),  ($\mathcal{F}$1) and ($\mathcal{F}$2)
 hold and there exists enough small $\delta\in (0,1)$ such that
 \begin{eqnarray}
   \label{dvv4}        F(t,x_1,x_2)
  &  =  &  \left(\left|\sin\frac{\pi}{4}t\right|+1\right)|x_1|^{\frac{3}{2}}+\left(\cos^2\frac{\pi}{4}t+1\right)|x_2|^{2}\nonumber\\
  & \ge &  |x_1|^{\frac{3}{2}}+|x_2|^{2} \nonumber\\
  & \ge &  |x_1|^{2}+|x_2|^{\frac{5}{2}}, \quad \forall \
  |x_1|<\delta, |x_2|<\delta.
 \end{eqnarray}
  Let $\beta_1=2$, $\beta_2=\frac{5}{2}$ and $M_1=M_2=1$. Then
  (\ref{dvv4}) implies that ($\mathcal{F}$3) holds. Hence,  by Theorem
  1.2, we obtained that system (\ref{eqq3}) has at least $48$ distinct pairs
  of 4-periodic
 solutions.

  \vskip4mm
  \noindent
 {\bf Acknowledgment}
  \par
  This project is supported by the National Natural
 Science Foundation of China (No: 11301235) and Tianyuan Fund for Mathematics of the National Natural
 Science Foundation of China (No: 11226135).

\textbf{}

 \end{document}